\documentclass{amsart}
\usepackage{amssymb, amscd}
\numberwithin{equation}{section}

\def\GG{{\mathbb G}}

\def\tto{\longrightarrow}

\def\Hcal{{\mathcal H}}

\def\Ocal{{\mathcal O}}
\def\Pcal{{\mathcal P}}

\newcommand\Hom{\operatorname{Hom}}
\newcommand\Ker{\operatorname{Ker}}

\newcommand\riso{\mathrel{\hskip2pt\raise-2.5pt\hbox{$\widetilde{\phantom{xx}}$}
\kern-16pt\longrightarrow}}

\newcommand\proofsquare{\nobreak\hfill \hbox{%
\vrule height 5pt
\kern-.4pt
 \vbox{%
\hrule width 5pt depth0pt height.4pt
 \kern4.6pt \hrule  }
\kern-3.75pt
\vrule height 5pt}\kern1pt
\par}

\newtheorem{theorem}{Theorem}[section]
\newtheorem{lemma}[theorem]{Lemma}
\newtheorem{proposition}[theorem]{Proposition}
\newtheorem{corollary}[theorem]{Corollary}

\newtheorem{definition-lemma}[theorem]{Definition-Lemma}

\theoremstyle{definition}
\newtheorem{definition}[theorem]{Definition}
\newtheorem{example}[theorem]{Example}

\theoremstyle{remark}

\def\llongrightarrow{\relbar\joinrel\relbar\joinrel\rightarrow}

\begin{document}

\title[
On the Height of Calabi-Yau Varieties in Positive Characteristic]
{ On the Height of
Calabi-Yau varieties \\
in Positive Characteristic}
\author{G.\ van der Geer}

\address{Faculteit Wiskunde en Informatica, University of
Amsterdam, Plantage Muidergracht 24, 1018 TV Amsterdam, The Netherlands}
\email{geer@science.uva.nl}
\author{T.\ Katsura}
\address{Graduate School of Mathematical Sciences\\
 University of Tokyo \\ Komaba, Meguro-ku\\
 Tokyo\\
153-8914 Japan}
\email{tkatsura@ms.u-tokyo.ac.jp}


\begin{abstract}
We study invariants of Calabi-Yau varieties in positive characteristic,
especially the height of the Artin-Mazur formal group. We illustrate
these results by Calabi-Yau varieties of Fermat and Kummer type. 
\end{abstract}

\maketitle
\begin{section}{Introduction}

The large measure of 
attention that complex Calabi-Yau varieties drew in recent years stands
in marked contrast to the limited attention for their counterparts in
positive characteristic. Nevertheless, we think these varieties
deserve a greater interest, especially since the special nature of these
varieties lends itself well for excursions into the largely 
unexplored territory of varieties in positive characteristic. 
In this paper we mean by a Calabi-Yau variety a smooth complete 
variety of dimension $n$ over a field with $\dim H^i(X,O_X)=0$ 
for $i=1,\ldots,n-1$ and with trivial canonical bundle. 
We study some invariants of Calabi-Yau varieties 
in characteristic $p>0$, especially the height $h$
of the Artin-Mazur formal group  for which we prove the estimate 
$h \leq h^{1,n-1} +1$ if $h\neq \infty$. 
We show how this invariant is related to
the cohomology of sheaves of closed forms. 

It is well-known that K3 surfaces do not possess non-zero global $1$-forms.
The analogous statement about the existence of global $i$-forms
with $i=1$ and $i=n-1$ on a $n$-dimensional Calabi-Yau variety
is not known and might well be false in positive characteristic. 
We show that for a Calabi-Yau
variety 
of dimension $\geq 3$
over an algebraically closed field $k$ of characteristic $p>0$ 
with no non-zero global $1$-forms there is no $p$-torsion
in the Picard variety and ${\rm Pic}/p{\rm Pic}$ is isomorphic to
$NS/pNS$ with $NS$ the N\'eron-Severi group of $X$. If in addition $X$ does not
have a non-zero global $2$-form then $NS/pNS\otimes_{{\bf F}_p} k$ maps injectively
into $H^1(X,\Omega_X^1)$. This yields the estimate $\rho \leq h^{1,1}$
for the Picard number. We also study Calabi-Yau varieties of Fermat type and
of Kummer type to illustrate the results.

\end{section}

\begin{section}{The Height of a Calabi-Yau Variety}
The most conspicuous invariant of a Calabi-Yau variety $X$ of dimension $n$
in characteristic $p>0$ is its {\sl height}. There are several ways to
define it, using crystalline cohomology or formal groups. In the latter
setting one considers the functor
$F^r_X: {\rm Art} \to {\rm Ab}$ 
defined on the category of local Artinian $k$-algebras 
with residue field $k$
by 
$$
F^r_X(S):= {\rm Ker} \{ H_{\rm et}^r(X \times S, \GG_m) 
\tto  H_{\rm et}^r(X,\GG_m) \}.
$$
According to a theorem of Artin and Mazur~\cite{AM}, for
a Calabi-Yau variety $X$ and $r=n$ this functor is representable by
a smooth formal group $\Phi_X$ 
of dimension~$1$ with tangent space $H^n(X,O_X)$.
Formal groups of dimension $1$ in characteristic $p>0$ 
are classified up to isomorphism
by their height $h$
which is a natural number $\geq 1$ or $\infty$. In the former case
($h\neq \infty$)
the formal group is $p$-divisible, while in the latter case the formal
group is isomorphic to the additive formal group $\hat{\GG}_a$.

For a non-singular complete variety $X$ over an algebraically closed field  $k$
of characteristic $p>0$  we let $W_mO_X$ be the coherent sheaf of 
Witt rings of length $m$. It has three operators $F$, $V$ and $R$
given by $F(a_0,\ldots,a_m)=(a_0^p,\ldots,a_m^p)$, $V(a_0,\ldots,a_m)=
(0,a_0,\ldots,a_{m})$ and $R(a_0,\ldots,a_m)=(a_0,\ldots,a_{m-1})$
satisfying the relations $RVF=FRV=RFV=p$. The cohomology groups
$H^i(X,W_mO_X)$ with the maps induced by $R$ form a projective system
of finitely generated $W_m(k)$-modules. The projective limit is the
cohomology group $H^i(X,WO_X)$. Note that this need not be a finitely generated
$W(k)$-module. It has semi-linear operators $F$ and $V$.

Let $X$ be a Calabi-Yau manifold of dimension $n$. The vanishing
of the groups $H^i(X,O_X)$ for $i\neq 0, n$ and the exact sequence 
$$
0 \to W_{m-1}O_X \to  W_mO_X \to O_X \to 0
$$ 
imply that $H^{i}(X,W_mO_X)$ vanishes for $i=1,\ldots, n-1$ and
all $m>0$, hence $H^{i}(X,WO_X)=0$. We also see that
restriction $R: W_mO_X \to W_{m-1}O_X$ 
induces a surjective map
$H^n(X,W_mO_X) \to H^n(X,W_{m-1}O_X)$ 
with kernel $H^n(X,O_X)$. The fact that $F$ and $R$ commute
implies that if the induced map 
$F: H^n(X,W_mO_X)\to H^n(X,,W_mO_X)$ vanishes then 
$F: H^n(X,W_iO_X)\to H^n(X,W_iO_X)$ vanishes for
$i<m$ too. It also follows that $H^n(X,W_iO_X)$ is a $k$-vector space
for $i<m$.

It is known by Artin-Mazur~\cite{AM} 
that the Dieudonn\'e module of the formal group $\Phi_X$ is
$H^n(X,WO_X)$ with $WO_X$ the sheaf of Witt vectors of $O_X$.
This implies the following result, cf.\ \cite{GK1} where 
we proved this
for K3-surfaces. We omit the proof which is similar to
that for K3 surfaces.

\begin{theorem}\label{th:characterization}
For a Calabi-Yau manifold $X$ of dimension $n$ we
have the following characterization of the height:
$$
h(\Phi_{X}) = \min \{ i \geq 1 \colon [F: H^{n}(W_iO_X) \rightarrow
H^{n}(W_{i}O_X)] \neq 0 \}.
$$
\end{theorem}

 We now connect this with de Rham cohomology.
Serre introduced in \cite{S} a map
$D_i: W_i(O_X) \rightarrow \Omega_X^1$ of sheaves in the following way:
$$
D_i(a_0,a_1,\ldots,a_{i-1}) =  a_0^{p^{i-1}-1}da_0 +\ldots +
a_{i-2}^{p-1}da_{i-2} + da_{i-1}.
$$
It satisfies $D_{i+1}V = D_{i}$,
and Serre showed that this induces an injective map of 
sheaves of additive groups
$$
D_i : W_iO_X/FW_iO_X \rightarrow \Omega^1_X \eqno(1)
$$
The exact sequence
$0 \rightarrow  W_iO_X \stackrel{F}{\longrightarrow}
W_iO_X \longrightarrow  W_iO_X/FW_ iO_X \rightarrow 0$
gives rise to an isomorphism
$$
H^{n - 1}(W_iO_X/FW_iO_X) \cong 
\Ker[F: H^{n}(W_iO_X) \rightarrow H^{n}(W_iO_X)].
\eqno(2)
$$
\begin{proposition}\label{prop:n-1-inj}
If $h\neq \infty$ then the induced map
$$D_{i} :H^{n-1}(X,W_iO_X/FW_iO_X) \to  H^{n-1}(X,\Omega^{1}_{X})$$
is injective, and $\dim {\rm Im} \, D_{i} =  \min \{i,h-1\}$.
\end{proposition}
\proof{We give a proof for the reader's convenience.
Take an affine open covering $\{ U_{i}\}$ of $X$. 
Assuming some $D_{\ell}$ is not injective, we let
$\ell$ be the smallest natural number such that
$D_{\ell}$ is not injective on $ H^{n-1}(W_{\ell}O_X/FW_{\ell}O_X)$.
Let $ \alpha = \{ f_I \}$ with $f_I=(f^{(0)}_I,\ldots , f^{(\ell -1)}_I)
 \in \Gamma(U_{i_0,\ldots,i_{n-1}},W_{\ell}O_X)$ represent
a non-zero element of of $H^{n-1}(W_{\ell}O_X/FW_{\ell}O_X)$ 
such that $D_{\ell}(\alpha)$ is zero in 
$H^{n-1}(\Omega^{1}_{X})$. Then there exists
elements $\omega_J=\omega_{j_{0}i_{1}\ldots j_{n-2}}$ in
$\Gamma (U_{j_{0}} \cap \ldots \cap U_{j_{n-2}}, \Omega_{X}^{1})$ such that
$$
\sum_{j=1}^{\ell} (f_I^{(j)})^{p^{\ell-j}} d\log f^{(j)}_I =
\sum_{j} \omega_{I_j},
$$
where the multi-index
$I_j=\{i_0,\ldots,i_{n-1}\}$ is obtained from $I$ by omitting $i_j$.
By applying the inverse Cartier operator we get an equation
$$
\sum_{j=1}^{\ell} (f_I^{(j)})^{p^{\ell-j+1}} d\log f^{(j)}_I  +
df_I^{(\ell)}= \sum_{j} \tilde{\omega}_{I_j}
$$
for certain functions $f_I^{(\ell)}$ and differential forms
$\tilde{\omega}_{I_j}$ with
$C(\tilde{\omega}_{I_j})=\omega_{I_j}$.
Since $\alpha$ is an non-zero element of $H^{n-1}(W_{\ell}O_X/FW_{\ell}O_X)$, 
the element $\beta =(f^{(0)}_I,\ldots , f^{(\ell -1)}_I, f^{(\ell)}_I)$
gives a non-zero element of $H^{n-1}(W_{\ell +1}O_X/FW_{\ell +1}O_X)$. 
In view of (2) for $i=\ell+1$ the element
$\beta$ gives a non-zero element $\tilde{\beta}$ of $H^{n}(W_{\ell +1}O_X)$
such that $F(\tilde{\beta}) = 0$
in $H^{n}(W_{\ell +1}O_X)$. 
Take the element $\tilde{\alpha}$ in $H^{n}(W_{\ell}O_X)$
which corresponds to the element $\alpha$ under the isomorphism
(2) for $i=\ell$. 
Then we have $F(\tilde{\alpha}) = 0$ in $H^{n}(W_{\ell}O_X)$, and
$R^{\ell}(\tilde{\alpha}) \neq 0$
in $H^{n}(X, {O}_{X})$ by the assumption on $\ell$. Therefore, we have
$R^{\ell}(\tilde{\beta}) \neq 0$ in $H^{n}(X, {O}_{X})$, and the
 elements 
$V^jR^j(\tilde{\beta})$ for $j=0, \ldots,\ell$ generate
$H^{n}(W_{\ell+ 1}O_X)$. 
Hence the Frobenius map is zero on $H^{n}(W_{\ell+ 1}O_X)$.
Repeating this argument, we conclude that the Frobenius map is zero on
$H^{n}(W_{i}O_X)$ for any $i>0$
and this contradicts the assumption $h \neq \infty$.}

\begin{corollary} If the height $h$ 
of an $n$-dimensional Calabi-Yau variety $X$
is not $\infty$ then 
$h \leq \dim H^{n-1}(\Omega^{1}_{X}) +1$.
\end{corollary}

\begin{definition}
A Calabi-Yau manifold $X$ is called {\sl rigid} 
if $\dim H^{n-1}(X,\Omega^{1}_{X})=0$.
\end{definition}

Please note that the tangent sheaf $\Theta_X$ is the dual of $\Omega_X^1$, 
hence by the triviality of the canonical bundle it is isomorphic 
to $\Omega_X^{n-1}$.
Therefore, by Serre duality the space of infinitesimal 
deformations $H^1(X,\Theta_X)$
is isomorphic to the dual of $H^{n-1}(X, \Omega_X^1)$.

\begin{corollary}
The height of a rigid Calabi-Yau manifold $X$
is either 1 or $\infty$.
\end{corollary}

\end{section}

\begin{section}{Cohomology Groups of Calabi-Yau Varieties}

Let $X$ be a Calabi-Yau variety of dimension $n$ over $k$. The existence
of Frobenius provides the de Rham cohomology with a very rich structure
from which we can read off characteristic $p$ properties. If $F: X \to X^{(p)}$
is the relative Frobenius operator then the Cartier operator $C$ gives
an isomorphism 
$$
{\Hcal}^j(F_*\Omega_{X/k}^{\bullet})=
\Omega_{X,\hbox{\rm $d$-closed}}^j/ d\Omega_X^{j-1}
\riso
\Omega^j_{X^{(p)}}$$
of sheaves on $X^{(p)}$. We generalize the sheaves $d\Omega_X^{j-1}$ and
$\Omega_{X,\hbox{\rm $d$-closed}}^j$ by setting (cf.\ \cite{I})
$$
B_0\Omega_X^j=(0), \quad B_1\Omega_X^j=d\Omega_X^{j-1}, 
\quad B_{m+1}\Omega_X^j= C^{-1}(B_m\Omega_X^j).
$$
and
$$
Z_0\Omega_x^j= \Omega_X^j, \quad 
Z_1\Omega_X^j=\Omega_{X,\hbox{\rm $d$-closed}}^j,
\quad Z_{m+1}\Omega_X^j= {\rm Ker}(dC^m).
$$
Note that we have the inclusions
$$
\begin{array}{c}
0=B_0 \Omega_X^j \subset B_1 \Omega_X^j \subset \ldots \subset B_m \Omega_X^j
 \subset \ldots \\
\ldots \subset Z_m \Omega_X^j
\subset \ldots \subset Z_{1} \Omega_X^j \subset Z_0 \Omega_X^j = \Omega_X^j
\end{array}
$$
and that we have an exact sequence
$$
0 \rightarrow Z_{m+1}\Omega_X^j \longrightarrow Z_m\Omega_X^j
\stackrel{dC^m}{\longrightarrow} d\Omega_X^j \rightarrow 0.
$$
Alternatively, 
the sheaves $B_m \Omega_X^j$ and $Z_m \Omega_X^j$ can be viewed as
locally free subsheaves of
$(F^m)_*\Omega_X^j$ on $X^{(p^m)}$.
Grothendieck duality implies that for every $j\geq 0$ there is a perfect
pairing  of $O_{X^{(p^m)}}$-modules $F_{*}^m\Omega_{X}^j 
\otimes F_*^m\Omega_X^{n-j}
\tto \Omega_{X^{(p^m)}}^n$ given by $(\alpha, \beta) \mapsto C^m(\alpha \wedge 
\beta)$. This induces perfect pairings of $O_{X^{(p^m)}}$-modules
$$
B_m\Omega_X^{j} \otimes F_*^m\Omega_X^{n-j} / Z_m \Omega_X^{n-j}
\to \Omega_{X^{(p^m)}}
$$
and
$$
Z_m\Omega^j_X \otimes F_*^m \Omega_X^{n-j}/B_m\Omega_X^{n-j}
\to \Omega_{X^{(p^m)}}
$$
Now we have an isomorphism $F_*^m \Omega^{j}_X/Z_m\Omega^{j}_X \cong
B_m\Omega^{j+1}_X$ induced by the map $d$. 
Going back to the interpretation of the $B_m\Omega_X^j$ as sheaves
on $X$ we find in this way for 
$1\leq j \leq n$ and $m>0$ perfect pairings
$$
B_m\Omega_X^j \otimes B_m \Omega_X^{n+1-j} \to \Omega_X^n
\qquad (\omega_1\otimes \omega_2) \mapsto C^m (\omega_1 \wedge \omega_2).
$$
We first note another interpretation for $B_m\Omega^1_X$:
the injective map of sheaves of additive groups
$D_m : W_m({\Ocal}_X)/FW_m({\Ocal}_X) \rightarrow \Omega^1_X$
induces an isomorphism
$$
D_m: W_m({\Ocal}_X)/FW_m({\Ocal}_X) \riso B_m\Omega_X^1. \eqno(3)
$$
We write $h^i(X,-)$ for $\dim_k H^i(X,-)$. Note that duality implies
$h^i(B_m\Omega_X^{n})=h^{n-i}(B_m\Omega_X^1)$.
\begin{proposition}\label{prop:Bn}
We have $h^i(X,B_m\Omega_X^1)=0$ unless $i=n$ or $i=n-1$.
If $i=n-1$ or $i=n$ we have
$$
h^{i}(B_m\Omega_X^1)  =
\begin{cases} \min\{ m, h-1\} & \hbox{\rm if $h \neq \infty $}\\
m & \hbox{\rm if $h=\infty$.} \\
\end{cases}
$$
\end{proposition}
\proof{The statement about $h^{n-1}(B_m\Omega_X^1)$ follows from (3)
and the characterization of the height given in Section~2.
The other statements follow from the long exact sequence associated 
with the short exact sequence
$$
0 \longrightarrow {\Ocal}_{X} \stackrel{F}{\longrightarrow} {\Ocal}_{X}
\stackrel{d}{\longrightarrow} d{\Ocal}_{X} \longrightarrow 0
$$
and the exact sequence 
$$
0 \to B_m \to B_{m+1} {\buildrel C^m \over \llongrightarrow} B_1 \to 0. \eqno(4)
$$
The details can safely be left to the reader.
This concludes the proof.}
\medskip

The natural inclusions
$B_i\Omega_{X}^{j} \hookrightarrow \Omega_{X}^{j}$
and 
$Z_i\Omega_{X}^{j} \hookrightarrow \Omega_{X}^{j}$
of sheaves of groups on $X$ induce homomorphisms
$$
H^1( B_i\Omega_{X}^{j}) \rightarrow H^1(\Omega_{X}^{j}) \quad
{\rm and} \quad H^1(Z_i\Omega_{X}^{j}) \rightarrow H^1(\Omega_{X}^{j})
$$
whose images are denoted by  ${\rm Im} \, H^1(B_i\Omega_{X}^{j})$ and
${\rm Im} \, H^1(Z_i\Omega_{X}^{j})$.
Note that we have a non-degenerate cup product pairing
$$
\langle \, , \, \rangle :  
H^{n-1}(X, \Omega_X^1) \otimes H^1(X, \Omega_X^{n-1}) \to H^n(X, \Omega_X^n)
\cong k.
$$
\begin{lemma}\label{lm:orthogonal}
The images 
${\rm Im}\,  H^{n-1}(B_i\Omega_{X}^{1})$ and ${\rm Im}\,  H^1(Z_i\Omega_{X}^{n-1})$ 
are orthogonal to each other for the pairing $\langle \, , \, \rangle$.
\end{lemma}
\proof{From the definitions it follows that for elements
$\alpha \in H^{n-1}(B_i\Omega_X^{1})$ and $\beta \in H^1(Z_i\Omega_X^{n-1})$
we have $C^i(\alpha \wedge \beta)=0$. The long exact sequence associated
to 
$$
0\to B_1\Omega_X^n \to Z_1\Omega_X^n \to \Omega_X^n \to 0
\eqno(5)
$$ 
together with the fact that $H^n(Z_i\Omega_X^n)=H^n(\Omega_X^n)$
for $i\geq 0$ implies that $C$ acts without kernel on $H^n(\Omega_X^n)$.
So the pairing is non-degenerate.
}

\begin{lemma}
If $h \neq \infty$ we have
$\dim {\rm Im} H^{1}(X, Z_{i}\Omega_{X}^{n - 1}) 
= \dim H^{1}(\Omega_{X}^{n-1}) - i$
for $0\leq i \leq h - 1$.
\end{lemma}
\proof{If the height $h = 1$ then we have
$H^{n-1}(B_i\Omega_{X}^{1}) = 0$ by (4) and moreover the vanishing of
$H^i(X,d\Omega_X^{n-1})$ and the exact sequence
$$
0 \rightarrow Z_{i+1}\Omega_{X}^{n-1} \longrightarrow Z_{i}\Omega_{X}^{n-1}
\stackrel{dC^{i}}{\llongrightarrow}
 d\Omega_{X}^{n-1}  \rightarrow 0 \eqno(6)
$$
imply that
${\rm Im} \, H^1(X,Z_i\Omega_X^{n-1})= H^1(X,\Omega_X^{n-1})$
for $i \ge 1$.
For $2 \leq h < \infty$, we know by Proposition~\ref{prop:n-1-inj}
that
${\rm Im}\,  H^{n-1}(X, B_{i}\Omega_{X}^{1})
\subset H^{n-1}(X, \Omega_{X}^{1})$
is of dimension $\min \{i, h - 1\}$.
The exact sequence (6) gives an exact sequence
$$
k \longrightarrow H^1(Z_{i+1}\Omega_{X}^{n-1})
\stackrel{\psi_{i+1}}{\llongrightarrow}
H^1(Z_{i}\Omega_{X}^{n-1}) \longrightarrow k
$$
from which we deduce that either
$\dim \psi_{i + 1}(H^1(Z_{i+1}\Omega_{X}^{n-1})) = 
\dim H^1(Z_{i}\Omega_{X}^{n-1}) + 1$
or
$\dim \psi_{i + 1}(H^1(Z_{i+1}\Omega_{X}^{n-1})) = 
\dim H^1(Z_{i}\Omega_{X}^{n-1})$.
By induction $\dim {\rm Im} H^1(Z_{i}\Omega_{X}^{n-1})$
is at least  $\dim H^1(\Omega_{X}^{n-1}) - i$. On the other hand, by 
Proposition~$\ref{prop:Bn}$
we have $\dim {\rm Im} H^1(Z_{i}\Omega_{X}^{n-1}) \leq
\dim H^1(\Omega_{X}^{n-1}) - i$ for $i \leq h - 1$.}

\begin{lemma}
If $X$ is a Calabi-Yau manifold of dimension $n$ with $h = \infty$ then
$$
({\rm Im}\, H^{n-1}(X, B_{i}\Omega_{X}^{1}))^{\perp} = 
{\rm Im}\, H^1(Z_{i}\Omega_{X}^{n-1}).
$$
\end{lemma}
\proof{We prove this by induction on $i$. 
By the exact sequence (5)
we have $\dim H^{i}(X, d\Omega_{X}^{n-1}) = 1$ for $i = 0, 1$.  Thus,
by the exact sequence (6) we see that the difference
$\dim {\rm Im} H^{1}(Z_{i}\Omega^{n-1}) - \dim {\rm Im} 
H^{1}(Z_{i+1}\Omega^{n-1})$
is equal to $0$ or $1$, and we have an exact sequence
$$
H^{1}(Z_{i+1}\Omega^{n-1}) \stackrel{\phi}{\longrightarrow} 
H^{1}(Z_{i}\Omega^{n-1}) \stackrel{dC^{i}}{\longrightarrow}
H^{1}(d\Omega_{X}^{n-1}).
$$
Assume that
${\rm Im}\, H^{n-1}(B_{j-1}\Omega^{1}) \neq {\rm Im}\, H^{n-1}(B_{j}\Omega^{1})$
for $j \leq i$ and
$
{\rm Im}\, H^{n-1}(B_{i}\Omega^{1})~=
{\rm Im}\, H^{n-1}(B_{i+1}\Omega^{1})
$.
By Lemma~\ref{lm:orthogonal},
$$
{\rm Im}H^{1}(Z_{i-1}\Omega^{n-1}) ~\supset~ 
{\rm Im}H^{1}(Z_{i}\Omega^{n-1})
$$
and ${\rm Im}H^{1}(Z_{i-1}\Omega^{n-1}) \neq {\rm Im}H^{1}(Z_{i}\Omega^{n-1})$ 
for $j \leq i$.
Suppose ${\rm Im}H^{1}(Z_{i}\Omega^{n-1}) \neq {\rm Im}
H^{1}(Z_{i+1}\Omega^{n-1})$.
The natural homomorphism $\phi : H^{1}(Z_{i+1}\Omega^{n-1}) \rightarrow
H^{1}(Z_{i}\Omega^{n-1})$ is not surjective. Since $H^{1}(d\Omega_{X}^{n-1}) 
\cong k$,
we see that $dC^{i} : H^{1}(Z_{i}\Omega^{n-1}) \rightarrow
H^{1}(d\Omega_{X}^{n-1})$ is surjective and we factor it as
$$
H^{1}(Z_{i}\Omega^{n-1}) \stackrel{C^{i}}{\longrightarrow} 
H^{1}(\Omega^{n-1})\stackrel{d}{\longrightarrow} H^{1}(d\Omega_{X}^{n-1}).
$$
Since $dC^{i}$ is surjective, $d$ is not the
zero map on $C^{i}(H^{1}(Z_{i}\Omega^{n-1}))$.
Therefore, we have
$$
C^{i}(H^{1}(Z_{i}\Omega^{n-1})) \not\subset {\rm Im}\, H^{1}(Z_{1}\Omega^{n-1}).
$$
Take an affine open covering of $X$, and take any \v{C}ech cocycle
$C^{i}(\eta) = \{C^{i}(\eta_{jk})\}$ of $C^{i}(H^{1}(Z_{i}\Omega^{n-1}))$
with respect to this affine open covering.
Take any element $\zeta \in H^{n-1}(B_{i}\Omega^{1})$. Then there exists
an element $\tilde{\zeta}$ such that $C^{i}(\tilde{\zeta})= \zeta$.
We consider the image of the element $C^{i}({\eta}) \wedge \zeta$
in $H^{n}(X, \Omega^{n}_{X})$. Then, we have
$$
C^{i}(\tilde{\eta}) \wedge \zeta  = 
C^{i}(\tilde{\eta}) \wedge C^{i}(\tilde{\zeta})
           = C^{i}({\eta} \wedge \tilde{\zeta})
$$
Since ${\rm Im} H^{n-1}(B_{2i}\Omega^{1}) 
= {\rm Im} H^{n-1}(B_{i}\Omega^{1})$,
the image of $\tilde{\zeta}$ in $H^{n-1}(\Omega^{1}_{X})$ is contained in
${\rm Im}H^{n-1}(B_{i}\Omega^{1})$. As ${\rm Im}
H^{1}(Z_{i}\Omega^{n-1}))$ is
orthogonal to ${\rm Im}H^{n-1}(B_{i}\Omega^{1})$, 
we see that $\eta \wedge \tilde{\zeta}$
is zero in  $H^{n}(X, \Omega^{n}_{X})$, and we have
$C^{i}({\eta} \wedge \tilde{\zeta}) = 0$ in $H^{n}(X, \Omega^{n}_{X})$.
Therefore, we see that the image of $C^{i}(H^{1}(Z_{i}\Omega^{n-1}))$ in
$H^{1}(\Omega^{n-1}_{X})$ is orthogonal to ${\rm Im}H^{n-1}(B_{i}\Omega^{1})$
and we have
$$
C^{i}(H^{1}(Z_{i}\Omega^{n-1}))  \subset~ {\rm Im}
H^{n-1}(B_{i}\Omega^{1})^{\perp} 
        \subset ~ {\rm Im} H^{1}(Z_{i}\Omega^{n-1}) 
     \subset~ {\rm Im} H^{1}(Z_{1}\Omega^{n-1}),
$$
a contradiction. Hence, we have
${\rm Im}H^{1}(Z_{i}\Omega^{n-1}) = {\rm Im}H^{1}(Z_{i+1}\Omega^{n-1})$.}

\smallskip

Collecting results we get the following theorem.

\begin{theorem} If $X$ is a Calabi-Yau variety of dimension $n$ and
height $h$ then for $i \leq h-1$ we have
$$
{\rm Im} H^{n-1}(X, B_{i}\Omega^{1}_X)^{\perp} 
= {\rm Im} H^1(X,Z_i\Omega_X^{n-1}).
$$
\par
\end{theorem}

One reason for our interest in the
spaces ${\rm Im}H^1(X,Z_i\Omega_X^{n-1})$ comes from the fact
that they play a role
as tangent spaces to strata in the moduli space as in the
analogous case of K3 surfaces, cf.\ \cite{GK1}. We intend to come back to
this in a later paper.
\end{section}

\begin{section}{Picard groups}
We suppose that $X$ is a Calabi-Yau variety of dimension $n\geq 3$.
We have the following result for the space of regular 1-forms.
\begin{proposition}\label{prop:H0(Zi)}
All global $1$-forms are indefinitely closed: for $i\geq 0$  we have
 $H^{0}(X,Z_{i}\Omega_{X}^{1}) = H^{0}(X,\Omega_{X}^{1})$.
The action of the Cartier operator on this space is semi-simple.
\end{proposition}
\proof{ Since the sheaves $B_i\Omega_X^1$ have non non-zero
cohomology in degree $0$ and $1$ the exact sequence
$$
0 \rightarrow B_{i}\Omega_{X}^{1} \longrightarrow
Z_{i}\Omega_{X}^{1}
\stackrel{C^{i}}{\longrightarrow}
 \Omega_{X}^{1}  \rightarrow 0.
$$
implies $\dim H^{0}(Z_{i}\Omega_{X}^{1}) = \dim H^{0}(\Omega_{X}^{1})$.
Since the natural map $H^{0}(Z_{i}\Omega_{X}^{1}) 
\rightarrow H^{0}(\Omega_{X}^{1})$
is injective,  we have $H^{0}(Z_{i}\Omega_{X}^{1}) = H^{0}(\Omega_{X}^{1})$.
The second assertion follows from $H^{0}(B_{i}\Omega_{X}^{1}) = 0$.}

It is well known that for 
a $p^{-1}$-linear semi-simple homomorphism  $\lambda$
on a finite-dimensional vector space $V$ the map $\lambda -{\rm id}_V$ 
is surjective.
This means that we have a basis of logarithmic differential forms $C\omega =
\omega$.
\begin{corollary}\label{cor:iso}
If ${\rm id}$ denotes the identity homomorphism on $H^{0}(X,\Omega_{X}^{1})$
the map $C - {\rm id} :H^{0}(X,\Omega_{X}^{1})
\rightarrow H^{0}(X,\Omega_{X}^{1})$ is surjective.
\end{corollary}

\begin{proposition}
Suppose that $X$ is a smooth complete variety for which all global $1$-forms
are closed and such that $C$ gives a bijection 
$H^0(X,Z_1\Omega_X^1) \longrightarrow H^0(X,\Omega_X^1)$.
Then we have an isomorphism
$$
H^0(X,\Omega_X^1) \cong {\rm Pic}(X)[p] \otimes_{\bf Z}k.
$$
\end{proposition}
\proof{ Let $L$ be a line bundle representing an element $[L]$ of order $p$
in ${\rm Pic}(X)$. Then there exists a rational function $g \in k(X)^*$
such that $(g)=pD$, where $D$ is a divisor corresponding to $L$. One 
observes now by a local calculation
that $dg/g$ is a regular $1$-form and thus defines an element
of $H^0(X,\Omega_X^1)$. Conversely, if $\omega$ is a global regular $1$-form with $C\omega = \omega$ then $\omega$ can be represented locally as $df_i/f_i$
with respect to some open cover $\{ U_i\}$. From the relation
$df_i/f_i=df_j/f_j$ we see $d\log(f_i/f_j)=0$ and this implies $d(f_i/f_j)=0$.
Hence we see that $f_i/f_j=\phi_{ij}^p$ form some $1$-cocycle $\{ \phi_{ij} \}$.
This cocycle defines a torsion element of order $p$ of ${\rm Pic}(X)$.}
These two maps are each others inverse and the result follows.

We denote by ${\rm Pic}(X)$ (resp.\ $NS(X)$) the Picard group (resp.\
 N\'eron-Severi group)
 of $X$. If $L$ is a line bundle 
with transition functions $\{ f_{ij}\}$ then $d\log f_{ij}$ represents
the first Chern class of $L$. In this way we can define a homomorphism
$$
\varphi_{1}: {\rm Pic}(X)\longrightarrow  H^{1}(Z_{1}\Omega_{X}^{1}), \quad
    [L]    \mapsto       c_1(L)= \{df_{ij}/f_{ij}\}
$$
which obviously factors through ${\rm Pic}/p{\rm Pic}$.

\begin{proposition}\label{prop:Pic}
The homomorphism
$\varphi_{1} :  {\rm Pic}(X)/p{\rm Pic}(X) 
\longrightarrow H^{1}(X,Z_{1}\Omega_{X}^{1})$ is injective.
\end{proposition}
\proof{We take an affine open covering $\{U_{i}\}$.
Suppose that there exists an element $[L]$ such that
$\varphi_{1} ([L]) = 0$. Then there exists a d-closed regular 1-form
$\omega_{i}$ on an affine open set $U_{i}$ such that
$df_{ij}/f_{ij} = \omega_{j} - \omega_{i}$ on $U_{i} \cap U_{j}$ and we have
$df_{ij}/f_{ij} = C(\omega_{j}) - C(\omega_{i})$. Therefore, we have
$\omega_{j}  - C(\omega_{j}) 
= \omega_{i} - C(\omega_{i})$ on $U_{i} \cap U_{j}$.
This shows that there exists an regular 1-form $\omega$ on $X$ such that
$\omega = \omega_{i} - C(\omega_{i})$ on $U_{i}$.
By Corollary~$\ref{cor:iso}$, there exists an element
$\omega^{\prime} \in H^{0}(\Omega_{X}^{1})$ such that
$(C - {\rm id})\omega^{\prime} = \omega$. 
Replacing $\omega_{i} + \omega^{\prime}$ by
$\omega_{i}$, we have
$$
 df_{ij}/f_{ij} = \omega_{j} - \omega_{i}
$$
with $C(\omega_{i}) = \omega_{i}$. Then, there exists an regular function
$f_{i}$ on $U_{i}$ such that $\omega_{i} = df_{i}/f_{i}$. So we have
$d\log f_{ij} = d \log (f_{j}/f_{i})$. 
Therefore, there exists a regular function
$\varphi_{ij}$ on $U_{i} \cap U_{j}$ such that $f_{ij} 
= (f_{j}/f_{i})\varphi_{ij}^{p}$. Thus $[L]$ is a $p$-th power.
We conclude that $\varphi_{1} :  {\rm Pic}(X)/p{\rm Pic}(X) 
\rightarrow H^{1}(Z_{1}\Omega_{X}^{1})$
is injective.}
\begin{proposition}
The natural homomorphism
$H^{1}(Z_{1}\Omega_{X}^{1}) \rightarrow H^{2}_{DR}(X)$ is injective.
\end{proposition}

\proof{Let $\{U_{i}\}$ be an affine open covering of $X$.
A \v{C}ech cocycle  $\{\omega_{ij}\}$ in $H^{1}(Z_{1}\Omega_{X}^{1})$
is mapped to $\{(0, \omega_{ij}, 0)\}$ in  $H^{2}_{DR}(X)$.
Suppose this element is zero in $H^{2}_{DR}(X)$. Then there exist  elements
$(\{f_{ij}\}, \{\omega_{i}\})$ with 
$f_{ij} \in \Gamma (U_{i} \cap U_{j},O_{X})$ and
$\omega_{i} \in \Gamma (U_{i}, \Omega_{X}^{1})$ such that
$$
f_{jk}-f_{ik} + f_{ij} = 0  \quad
  \omega_{ij} = df_{ij} + \omega_{j} - \omega_{i} \quad
         d\omega_{i} = 0.
$$
Since $\{f_{ij}\}$ gives an element of $H^{1}(O_{X})$ and 
$H^{1}(O_{X}) = 0$, there exists an element $\{f_{i}\}$ 
such that $f_{ij} =f_{j} - f_{i}$ on $U_{i} \cap U_{j}$. 
Therefore, we have
$ \omega_{ij} = (df_{j} + \omega_{j}) - (df_{i} + \omega_{i})$.
Since $d(df_{i} + \omega_{i}) = 0$, we conclude that  
$\{\omega_{ij}\}$ is zero in $H^{1}(Z_{1}\Omega_{X}^{1})$.}

\smallskip

The results above imply the following theorem.

\begin{theorem}
The natural homomorphism
${\rm Pic}(X)/p{\rm Pic}(X) \longrightarrow H^{2}_{DR}(X)$ is injective.
\end{theorem}

\begin{lemma}
For a Calabi-Yau manifold $X$ of dimension $n \geq 3$ with non non-zero global
$1$-forms 
${\rm Pic}(X)$ has no $p$-torsion. Moreover,
$NS(X)$ has no $p$-torsion and
$$
  {\rm Pic}(X)/p {\rm Pic} (X) \cong NS (X)/ p NS (X).
$$
\end{lemma}
\proof{Take an affine open covering $\{U_{i}\}$ of $X$. 
Assume $\{ f_{ij} \}$ represents an element $[L] \in {\rm Pic}(X)$ which
is  $p$-torsion. Then, there exist regular functions 
$f_{i} \in H^{0}(U_{i}, O_{X}^{*})$ such that
$f_{ij}^{p} = f_{i}/f_{j}$.
The $df_{i}/f_{i}$ on $U_{i}$ glue together to yield a 
regular 1-form $\omega$ on $X$.
Since $H^{0}(X, \Omega_{X}^{1}) = 0$, we see $\omega = 0$, i.e., $df_{i} = 0$.
Therefore, there exist regular functions $g_{i} \in H^{0}(U_{i}, O_{X}^{*})$
such that $f_{i} = g_{i}^{p}$. Hence, we have $\{ f_{ij} \} \sim 0$ and we see
that ${\rm Pic} (X)$ has no $p$-torsion.
Since $H^{1}(X, O_{X}) = 0$, the 
group  scheme ${\rm Pic}(X)$ is reduced and 
${\rm Pic}^{0} (X) = 0$.
Since $NS(X)={\rm Pic}(X)/{\rm Pic}^{0}(X) \cong {\rm Pic}(X)$, 
the lemma follows.
}

\begin{lemma}\label{lm:NS}
Let $X$ be a Calabi-Yau manifold $X$ of dimension $n \geq 3$ with 
no non-zero global $2$-forms. Then, the homomorphism
$$
NS(X)/p NS(X) \cong {\rm Pic}(X)/p {\rm Pic}(X) 
\longrightarrow H^{1}(\Omega_{X}^{1})
$$
defined by $\{ f_{ij} \} \mapsto \{df_{ij}/f_{ij}\}$ is injective.
\end{lemma}
\proof{By the assumption $H^{0}(X, \Omega_{X}^{2}) = 0$
we have $H^{0}(X, d\Omega_{X}^{1}) = 0$. Therefore, from the exact sequence
$$
 0 \rightarrow Z_{1}\Omega_{X}^{1} \longrightarrow
 \Omega_{X}^{1} \stackrel{d}{\longrightarrow} d\Omega_{X}^{1} \rightarrow 0,
$$
we deduce a natural injection
$H^{1}(Z_{1}\Omega_{X}^{1}) \longrightarrow H^{1}(\Omega_{X}^{1})$.
So the result follows from Lemma~\ref{prop:Pic}.
}

\begin{theorem}
Let $X$ be a Calabi-Yau manifold $X$ of dimension $n \geq 3$ with
$H^{0}(X, \Omega_{X}^{i}) = 0$  for $i = 1, 2$. Then the natural homomorphism
$$
NS(X)/p NS(X) \otimes_{{\bf F}_{p}} k \longrightarrow H^{1}(\Omega_{X}^{1})
\subset H^2_{dR}(X)
$$
is injective and  the Picard number satisfies
$\rho \leq \dim_{k} H^{1}(\Omega_{X}^{1})$.
\end{theorem}
\proof{Suppose that this homomorphism is not injective.
Then with respect to a suitable affine open covering $\{U_{i}\}$.
there exist elements $\{f_{ij}^{(\nu)}\}$ representing
non-zero elements in $NS(X)/p NS(X)$,
such that
$$
\sum_{\nu=1}^{\ell} a_{\nu}df_{ij}^{(\nu)}/f_{ij}^{(\nu)} 
= 0 \quad \hbox{in $H^1(\Omega_{X}^{1})$}
$$
for suitable $a_{\nu} \in k$.
We take such elements with the minimal $\ell$. 
We may assume $a_{1} = 1$
and we have
$a_{i}/a_{j} \notin {\bf F}_{p}$ for $i \neq j$. By Lemma~\ref{lm:NS}, we have
$\ell \geq 2$.
There exists
$\omega_{i} \in H^{1}(U_{i}, \Omega_{X}^{1})$ such that
$$
\sum_{\nu=1}^{\ell} a_{\nu}df_{ij}^{(\nu)}/f_{ij}^{(\nu)} 
= \omega_{j} -\omega_{i} \eqno(1)
$$
on $U_{i} \cap U_{j}$. There exists an element
${\tilde{\omega}}_{i} \in H^{1}(U_{i}, \Omega_{X}^{1})$ such that
$C({\tilde{\omega}}_{i}) = \omega_{i}$. Therefore, taking the Cartier inverse,
we have
$$
   \sum_{\nu=1}^{\ell} {\tilde{a}}_{\nu}df_{ij}^{(\nu)}/f_{ij}^{(\nu)} 
+ dg_{ij} =
{\tilde{\omega}}_{j} -{\tilde{\omega}}_{i}
$$
with ${\tilde{a}}_{\nu} \in k$, ${\tilde{a}}_{\nu}^{p} = a_{\nu}$, and
suitable $dg_{ij} \in H^{0}(U_{i}\cap U_{j}, dO_{X})$.
Since $\{ df_{ij}^{(\nu)}/f_{ij}^{(\nu)}\}$ is a cocycle, we see that
$\{dg_{ij}\} \in H^{1}(X, dO_{X})$. 
Since $H^{1}(X, dO_{X}) = 0$, there exists
an element $dg_{i} \in H^{0}(U_{i}, dO_{X})$ such 
that $dg_{ij} = dg_{j} - dg_{i}$.
Therefore, we have
$$
\sum_{\nu=1}^{\ell} {\tilde{a}}_{\nu}df_{ij}^{(\nu)}/f_{ij}^{(\nu)} =
({\tilde{\omega}}_{j} -dg_{j}) - ({\tilde{\omega}}_{i} -dg_i). \eqno(2)
$$
Subtracting (2) from (1) we get a non-trivial 
linear relation with a smaller $\ell$ in
$H^{1}(\Omega_{X}^{1})$, a contradiction.
}

\noindent
{\bf Remark.} In the case of a K3 surface $X$ the natural homomorphism
$$
NS(X)/p NS(X) \otimes_{{\bf F}_{p}} k \longrightarrow H_{DR}^{2}(X)
$$
is not injective if $X$ is supersingular in the sense of Shioda. Ogus
showed that the kernel can be used for describing the moduli
of supersingular K3 surfaces, cf. Ogus\cite{Og2}. 
So the situation is completely different in dimension $\geq 3$.
\medskip
\end{section}

\begin{section}{Fermat Calabi-Yau manifolds}

Again $p$ is a prime number and $m$ a positive integer which is prime to $p$.
Let $f$ be a smallest power of $p$ such that $p^{f} \equiv 1~\mbox{mod}~m$
and put $q = p^{f}$. We denote by ${\bf F}_q$ a  finite field of 
cardinality $q$.
We consider the Fermat variety $X_{m}^{r}(p)$ over ${\bf F}_{q}$ defined by
$$
 X_{0}^{m} + X_{1}^{m} + \ldots + X_{r + 1}^{m} = 0
$$
in projective space ${\bf P}^{r + 1}$ of dimension $r + 1$.
The zeta function of $X_{m}^{r}$ over ${\bf F}_{q}$ was calculated
by A. Weil (cf. \cite{We}). The result is:
$$
Z(X_{m}^{r}/{\bf F}_{q}, T)  = 
\frac{ P(T)^{(-1)^{r-1}}}{(1-T)(1-qT) \dots (1 - q^{r}T)},
$$
where $P(T) = \prod_{\alpha} (1- j(\alpha)T)$ with the product taken over
a set of vectors $\alpha$ and $j(\alpha)$
is a Jacobi sum defined as follows.
Consider the set
$$
A_{m, r} = \{ (a_{0}, a_{1},\ldots,a_{r+1}) \in {\bf Z}^{r+2} 
~\mid ~0< a_{i}  < m, {\textstyle \sum_{j=0}^{r+1}} a_{j} \equiv 0 (\bmod \,
m) \},
$$
and choose a character $\chi :{\bf F}_{q}^{*} \rightarrow {\bf C}^{*}$
of order $m$. For $\alpha =(a_0,a_1,\ldots,a_{r+1}) \in A_{m,r}$ we define
$$
j(\alpha) = (-1)^{r}\sum \chi(v_{1}^{a_{1}})\ldots \chi(v_{r+1}^{a_{r+1}}),
$$
where the summation runs over $v_{i} \in {\bf F}_{q}^*$ with
$1 + v_{1} + \ldots + v_{r+1} = 0$. Thus the 
$j(\alpha)$'s are eigenvalues of the Frobenius
map over ${\bf F}_{q}$ on the ${\ell}$-adic \'{e}tale cohomology group
$H_{et}^{r}(X_m^r, {\bf Q}_{\ell})$.

Now, let $\zeta = \exp(2\pi i/m)$ be a primitive $m$-th root of unity, 
and $K = {\bf Q}(\zeta)$ the corresponding cyclotomic field 
with Galois group $G= {\rm Gal}(K/{\bf Q})$.
For an element $t \in ({\bf Z}/m{\bf Z})^{*}$
we let $\sigma_{t}$ be the automorphism of $K$  defined by 
$\zeta \mapsto \zeta^{t}$.
The correspondence $t \leftrightarrow \sigma_t$ defines
an isomorphism $({\bf Z}/m{\bf Z})^{*} \cong G$ 
and we shall identify $G$ with $({\bf Z}/m{\bf Z})^{*}$ by this isomorphism.
We define a subgroup $H$ of order $f$ of $G$ by 
$H = \{p^{j}~\mbox{mod}~ m ~\mid~ 0 \leq j < f\}$

Let $\{ t_{1}, \ldots, t_{g}\}$ with $t_i \in {\bf Z}/m{\bf Z}^*$,
be a complete system of representatives of 
$G/H$ with $g = \vert G/H \vert$,
and put
$$
A_{H}(\alpha) = \sum_{t\in H}[\sum_{j=1}^{r+1}\langle ta_{j}/m\rangle],
$$
where $[a]$ (resp.\ $\langle a \rangle$) means the integral part 
(resp.\ the fractional part) of a rational number $a$. 

Choose a  prime ideal $\Pcal$ in K lying over $p$; 
it has  norm $N({\Pcal})=p^{f}= q$.
If ${\Pcal}_{i}$ denotes the prime ideal
${\Pcal}^{\sigma^{-1}_{-t_{i}}}$ we have the prime decomposition
$
(p) = {\Pcal}_{1}\cdots {\Pcal}_{g}
$
in $K$ and Stickelberger's theorem tells us that
$$
(j(\alpha)) = \prod_{i = 1}^{g} {\Pcal}_{i}^{A_{H}(t_{i}\alpha)},
$$
where $t_{i}\alpha = (t_{i}a_{0}, \ldots, t_{i}a_{r+1})$.
For the details we refer to
 Lang\cite{L} or Shioda-Katsura\cite{SK}.

Now we restrict our attention
to Fermat Calabi-Yau manifolds $X_{m}^{r}(p)$ with $m = r + 2$.

\begin{theorem}
Assume $r \geq 2$. Let $\Phi^{r}$ be the Artin-Mazur formal group of
 the $r$-dimensional Calabi-Yau variety
$X = X_{r+2}^{r}(p)$. The height $h$
of $\Phi^{r}$ is equal to either $1$ or $\infty$. Moreover, $h = 1$ if and only
if $p \equiv 1 ~{(\bmod}~r+2)$.
\end{theorem}

Before we give the proof of this theorem we state a technical lemma.

\begin{lemma}\label{lm:0}
Under the notation above, assume $[\sum_{j= 1}^{r+1}\langle ta_{j}/(r+2)
\rangle ] = 0$
with
$t \in ({\bf Z}/(r+2){\bf Z})^{*}$.
Then $a_{j} = t^{-1}$ in $({\bf Z}/(r+2){\bf Z})^{*}$ for all $j = 0, 1,
\ldots, r+1$.
\end{lemma}
\proof{Since $t \in ({\bf Z}/(r+2){\bf Z})^{*}$, we have $\langle
ta_{j}/(r+2)\rangle \geq
 1/(r+2)$.
Suppose there exists an index 
$i$ such that $ta_{i} \not\equiv 1 ~(\bmod\, r+2)$. Then we have the inequality
$\langle ta_{i}/(r+2)\rangle  \geq 2/(r+2)$ and thus
$\sum_{j= 1}^{r+1}\langle ta_{\rho}/(r+2)\rangle  \geq 1$, 
which contradicts the assumption.
So we have $ta_{i} \equiv 1 ~(\bmod \, r+2)$ and
$a_{j} \equiv t^{-1}$ for $j = 1, \ldots, r+1$.
Since $a_{0} + a_{1} + \ldots +a_{r+1} \equiv 0~ (\bmod \, r+2)$,
we conclude $a_{0} \equiv t^{-1}$.
}

{\sl Proof of the theorem.}
The Dieudonn\'{e} module $D(\Phi^{r})$ of $\Phi^{r}$
is isomorphic to $H^{r}(X, WO_X)$. 
We denote by Q(W) the quotient field
of the Witt ring $W(k)$ of $k$. Then, if $h < \infty$,  we have
$$
h = \dim_{Q(W)} H^{r}(X, WO_X) \otimes_{W(k)} Q(W).
$$ 
and by Illusie~\cite{I} we know
we have
$$
H^{r}(X, WO_{X}) \otimes_{W(k)} Q(W) 
\cong H^{r}_{\rm cris}(X)\otimes Q(W)_{[0, 1[}.
$$
According to
 Artin-Mazur~\cite{AM}, the slopes of $H^{r}_{\rm cris}(X)\otimes Q(W)$
are given by $({\rm ord}_{\Pcal} q)/f$ and  the
$({\rm ord}_{\Pcal} j(\alpha))/f$.
Hence, the height $h$ is equal to the number of $j(\alpha)$ 
such that $A_{H}(\alpha) < f$.

First, assume $p \equiv 1~(\bmod \, r+2)$, i.e.\ $f=1$. Then  
$H = \langle 1 \rangle$ and
$A_{H}(\alpha) < f = 1$ implies $A_{H}(\alpha)=0$. 
Therefore, by Lemma~\ref{lm:0},
we have $a_{j} = 1$ for all $j= 0, 1, \ldots, r+1$ and there is
only one $\alpha$, namely $\alpha = (1, 1, \ldots, 1)$, such that
${\rm ord}_{\Pcal}j(\alpha) = 0$. So we conclude $h = 1$ in this case.

Secondly, assume $p \not\equiv 1~(\bmod \, r+2)$. 
By definition, we have $f \geq 2$.
We now prove that there exists no $\alpha$ such that $A_{H}(\alpha)<f$.
Suppose $A_{H}(\alpha) = \sum_{t\in H}[\sum_{j= 1}^{r+1}\langle
ta_j/(r+2)\rangle ] < f$.
Then there exist an element $t \in H$ such that
$[\sum_{j= 1}^{r+1}\langle ta_{j}/(r+2)\rangle ] = 0$. 
By Lemma~\ref{lm:0} we have
$\alpha = (t^{-1}, t^{-1}, \ldots, t^{-1})$. For $t^{\prime} \in H$ with
$t^{\prime} \neq t$ we have 
$[\sum_{\rho = 1}^{r+1}\langle t^{\prime}t^{-1}/(r+2)\rangle] \neq 0$.
Therefore the inequality yields
$[\sum_{j= 1}^{r+1}\langle t^{\prime}t^{-1}/(r+2)\rangle]
 =1$ for $t^{\prime} \in H$,
$t^{\prime} \neq t$. Since $A_{H}(t\alpha) = A_{H}(\alpha)$
for any $t \in H$, by a translation by $t$, we may 
assume $\alpha = (1, 1, \ldots, 1)$, i.e., $t = 1$.
Moreover, we can take a representative of $t^{\prime} \in H$ such that
$0 < t^{\prime} < r+2$. Then,
\begin{align} \notag
1 &= [\sum_{j= 1}^{r+1}\langle t^{\prime}t^{-1}/(r+2)\rangle] 
   = [\sum_{j = 1}^{r+1}\langle t^{\prime}/(r+2)\rangle ]
 = [\sum_{j= 1}^{r+1}t^{\prime}/(r+2)]  \\
 &  = [(r+1)t^{\prime}/(r+2)] \notag
\end{align}
and we get $1\leq (r+1)t^{\prime}/(r+2) < 2$. By this inequality,
we see $t^{\prime} = 2$. 
Therefore, we have $H = \{1, 2\}$. Since $H$ is a subgroup
of $({\bf Z}/(r+2){\bf Z})^{*}$, we see that $2^{2} \equiv 1~ \mbox{mod}~r+2$.
Therefore, we have $r =1$, which contradicts our assumption.

Hence there exists no $\alpha$ such that 
${\rm ord}_{\Pcal} j(\alpha) < 1$ and 
we conclude $h = \infty$ in this case. This completes the proof
of the theorem.

\smallskip

For K3 surfaces we have two notions of supersingularity.
We generalize these to higher dimensions.

\begin{definition}
A Calabi-Yau manifold $X$ of dimension $r$ is said to be 
of {\sl additive Artin-Mazur type} (`supersingular in  the sense of Artin') 
if the height of Artin-Mazur formal group associated
with $H^{r}(X,{O}_{X})$ is equal to $\infty$.
\end{definition}

\begin{definition}
A non-singular complete algebraic variety $X$ of dimension $r$ is said to be
{\sl fully rigged} (`supersingular in the sense of Shioda')  
if all the even degree
{\'e}tale cohomology groups are spanned by algebraic cycles.
\end{definition}

By the theorem above, we know that the Fermat Calabi-Yau manifolds 
are of additive Artin-Mazur type if
and only if $p \not\equiv 1~\mbox{mod} ~m$ with $m = r +2$.
As to being fully rigged we have the following theorem.

\begin{theorem}[Shioda-Katsura~\cite{SK}]
Assume $m \geq 4$, $(p, m) = 1$ and $r$ is even. Then 
the Fermat variety $X_{m} ^{r}(p)$
is fully rigged if and only if there exists a positive integer $\nu$ such that
$p^{\nu} \equiv -1 ~{\mbox{mod}}~ m$.
\end{theorem}

M.\ Artin conjectured that a K3 surface $X$ is supersingular 
in the sense of Artin if and only if $X$ is supersingular in the sense 
of Shioda. He also showed that ``if part" holds. 
In the case of the Fermat K3 surface, i.e, $X_{4}^{2}(p)$, by the two
theorems above, we see, as is well-known, that the Artin conjecture holds.
However, in the case of even $r \geq 4$, the above two theorems 
imply that this straightforward
generalization of the Artin conjecture to higher 
dimension does not hold.
\end{section}

\begin{section}{Kummer Calabi-Yau manifolds}
Let $A$ be an abelian variety of dimension $n \geq 2$ defined over 
an algebraically closed field of characteristic $p >0$, 
and $G$ be a finite group
which acts on $A$ faithfully. Assume that the order of $G$ is prime to $p$,
and that the quotient variety $A/G$ has a resolution which is 
a Calabi-Yau manifold $X$. We call $X$ a Kummer Calabi-Yau manifold.
We denote by $\pi : A \longrightarrow A/G$
the projection, and by $\nu: X \longrightarrow A/G$ the resolution.

\begin{theorem}
Under the assumptions above the Artin-Mazur
formal group $\Phi_{X}^{n}$ is isomorphic to the Artin-Mazur
formal group $\Phi_{A}^{n}$.
\end{theorem}
\proof{Since the order of $G$ is prime to $p$, the singularities of $A/G$
are rational, and we have $R^{i}{\nu}_{*}{O}_{X} = 0 $ for $i \geq 1$.
So by the Leray spectral sequence we have 
$H^{n}(A/G, {O}_{A/G})\cong H^{n}(X,O_{X}) \cong k$ and
$H^{n-1}(A/G,{O}_{A/G})\cong H^{n-1}(X,O_{X}) \cong 0$.
It follows that the Artin-Mazur formal group 
$\Phi_{A/G}^{n}$ is pro-representable
by a formal Lie group of dimension $1$ (cf.\  Artin-Mazur\cite{AM}). 
Since the tangent space $H^{n}(A/G,O_{A/G})$ of $\Phi_{A/G}^{n}$
is naturally isomorphic to the tangent space 
$H^{n}(X, {O}_{X})$ of $\Phi_{X}^{n}$ 
as above, the natural homomorphism from $\Phi_{A/G}^{n}$ to $\Phi_{X}^{n}$
is non-trivial. One-dimensional formal groups are classified by their 
height and between formal groups of different height 
there are no non-trivial homomorphisms .
So the height of $\Phi_{A/G}^{n}$ is equal to that of $\Phi_{X}^{n}$ 
and we thus see that $\Phi_{A/G}^{n}$ and $\Phi_{X}^{n}$ are isomorphic. 

Since the order of $G$ is prime to $p$, there is a non-trivial trace map
from $H^{n}(A, {O}_{A})$ to $H^{n}(A/G, {O}_{A/G})$. Therefore,
$\pi^{*}: H^{n}(A/G, {O}_{A/G}) \longrightarrow H^{n}(A, {O}_{A})$
is an isomorphism. Therefore, as above we see
that the height of $\Phi_{A/G}^{n}$ is equal to the height of $\Phi_{A}^{n}$,
and that $\Phi_{X}^{n}$ is isomorphic 
to $\Phi_{A}^{n}$.} Q.e.d.
\smallskip

Though the following lemma might be well-known to specialists
we give here a proof for the reader's convenience.

\begin{lemma}\label{lm:abelian} Let $A$ be an abelian variety
of dimension $n \geq 2$ and $p$-rank $f(A)$.
The height $h$ of the Artin-Mazur formal group 
$\Phi_{A}$ of $A$ is as follows:
\begin{itemize}
\item[$({\rm 1})$] $h = 1$ if $A$ is ordinary, i.e., $f(A)=n$,
\item[$({\rm 2})$] $h = 2$ if $f(A) = n-1$,
\item[$({\rm 3})$] $h = \infty$ if $f(A) \leq n-2$.
\end{itemize}
\end{lemma}
\proof{We denote by $H^{i}_{\rm cris}(A)$ the $i$-th cristalline cohomology of
$A$ and as usual by $H^{i}_{\rm cris}(A)_{[\ell, \ell  + 1[}$ the
additive group of elements in $H^{i}_{\rm cris}(A)$ whose slopes are in the
interval $[\ell,\ell + 1[$. By the general theory
in Illusie \cite{I}, we have 
$$
H^{n}(A,W(O_A)) \otimes_{W}Q(W) \cong 
(H^{n}_{\rm cris}(A) \otimes_{W}Q(W))_{[0, 1[}
$$
with $Q(W)$ the quotient field of $W$. The theory of Dieudonn\'e modules
implies 
$$
 h = \dim_{Q(W)} D(\Phi_{A}) = \dim_{Q(W)} H^{n}(A,W(O_A)) \otimes_{W}Q(W)
 \quad 
\hbox{if $h < \infty$},
$$
and $\dim_{Q(W)} D(\Phi_{A}) = 0$ if $h = \infty$.
We know the slopes of $H^{1}_{\rm cris}(A)$ for each case. Since
we have 
$$
  H^{n}_{\rm cris}(A) \cong \wedge^{n} H^{1}_{\rm cris}(A),
$$
counting the number of slopes in $[0, 1[$ of 
$H^{n}_{\rm cris}(A)$ gives the result.}

\begin{corollary}
Let $X$ be a Kummer Calabi-Yau manifold of dimension $n$ obtained 
from an abelian variety $A$ as above.
Then the height of the Artin-Mazur formal group $\Phi_{X}^{n}$ is equal to 
either 1, 2 or $\infty$.
\end{corollary}

\begin{example}
Assume $p \geq 3$. Let $A$ be an abelian surface and $\iota$ the 
map $A \to A$ sending $a \in A$ to its inverse $-a\in A$.
We denote by $Km(A)$ the Kummer surface of A, i.e., the minimal resolution of 
$A/\langle \iota \rangle $.
Then $\Phi_{Km(A)}^{2}$ is isomorphic to $\Phi_{A}^{2}$.
\end{example}

\begin{example}\label{ex;abelian}
Assume $p \geq 5$, and let $\omega$ be a primitive third root of unity.
Let $E$ be a non-singular complete model of the elliptic curve
defined by $y^{2} = x^{3} + 1$, and let $\sigma$ be an automorphism of $E$
defined by $x\mapsto \omega x, ~y \mapsto y$. We set $A = E^{3}$ and 
put $\tilde{\sigma} = \sigma \times \sigma \times \sigma$. 
The minimal resolution $X$ of $A/\langle \tilde{\sigma} 
\rangle$ is a Calabi-Yau manifold, and 
the Artin-Mazur formal group $\Phi_{X}^{3}$ is isomorphic to $\Phi_{A}^{3}$.
\end{example}

Let $\omega$ be a complex number with positive imaginary part,
and $L = {\bf Z} + {\bf Z}\omega$ be a lattice in the complex
numbers ${\bf C}$.
From here on, we consider an elliptic curve $E = {\bf C}/L$, and we assume 
that $E$ has a model defined over an algebraic number field $K$. 
Then  $A = E \times E \times E$ is an abelian threefold, 
and we let $G\subseteq {\rm Aut}_K(E)$ 
be a finite group which faithfully acts on $A$. We assume that $G$ has only
isolated fixed points on $A$ and that the quotient variety $A/G$ has a crepant
resolution $\nu: X \to A/G$ defined over $K$, \cite{Reid}.
We denote by $\pi$ the projection 
$A \rightarrow A/G$.
For a prime $p$ of $K$, we denote by ${\bar{X}}$ the
reduction modulo $p$ of $X$.
 
\begin{example}\label{ex;omega}$[$K. Ueno, \cite{Ueno}$]$
Assume that $E$ is an elliptic curve defined over ${\bf Q}$
having complex multiplication $\sigma: E \to E$
by a primitive third root of unity.
Then $G={\bf Z}/3{\bf Z}=\langle \sigma \rangle$
acts diagonally on $A=E^3$.
A crepant resolution of $A/G$
gives a rigid Calabi-Yau manifold defined over ${\bf Q}$. 
For a prime number $p\geq 5$
the reduction modulo $p$ of $A$ is the abelian threefold given
in Example~\ref{ex;abelian}.
\end{example}

\begin{theorem}\label{intermediatejac}
Let $X$ be a Calabi-Yau obtained as crepant resolution of $A/G$
as above. Assume moreover that $X$ is rigid.
Then the elliptic  curve $E$ has complex multiplication
 and the intermediate Jacobian of $X$ is isogenous to $E$.
Moreover, if $\omega$ is a quadratic integer,
then the intermediate Jacobian of X
is isomorphic to $E$.
\end{theorem}

\begin{corollary}
Under the assumptions as in the theorem, we take a prime $p$ of good
reduction for $X$and let $\bar{X}$ be the reduction of $X$ modulo $p$. 
Then the height of the formal group $\Phi_{{\bar X}}$ is either $1$ or
$\infty$. It is $\infty$
if and only if the reduction of the intermediate Jacobian variety of $X$ 
at $p$ is a supersingular elliptic curve.
\end{corollary}

\begin{example}
We consider the reduction $\bar{X}$ modulo $p$ of the variety $X$
in the Example~\ref{ex;omega}.
We assume the characteristic of the residue field of $p$ 
is not equal to $2$ and $3$. 
Then the height $h(\Phi_{\bar{X}}) = \infty$
if and only if the reduction modulo $p$ of
the intermediate Jacobian of $X$ is a supersingular elliptic curve, and
this is the case if and only if $p \equiv 2~(\bmod \, 3)$.
\end{example}

Before we prove the theorem we introduce some notation.
We have a natural identification $H_{1}(E, {\bf Z}) = {\bf Z} + {\bf Z}\omega$.
Fixing a non-zero regular differential form $\eta$ on $E$ determines
a regular three form $\eta \times \eta \times \eta =\Omega_{A}$ on $A$.
We have a natural  homomorphism
$$
H^{3}(X, {\bf Z}) \rightarrow H_{dR}^{3}(X) \rightarrow H^{3}(X, O_{X}).
$$
If $X$ is rigid, the corresponding quotient $H^{3}(X,O_{X})/H^{3}(X,{\bf Z})$ 
gives the intermediate Jacobian of $X$. Since $\dim_{\bf C} H^{3}(X,O_{X}) = 1$,
the intermediate Jacobian of $X$ is isomorphic to an elliptic curve.

We can define the period map with respect to $\Omega_{A}$: 
$$
\pi_{A} : H_{3}(X, {\bf Z}) \longrightarrow  {\bf C}\qquad
           \gamma  \mapsto  \int_{\gamma} \Omega_{A}.
$$
By Poincar\'e duality we can identify $\pi_A$  with the natural projection
$H^3(X,{\bf Z}) \to H^3(X,{\bf C}) \to H^3(X,O_X)={\bf C}$
(cf.\ Shioda \cite{Sh1}, for instance.) 
There exists a regular 3-form $\Omega_{X}$ on ${X}$ such that 
$\Omega_{A} = (\nu^{-1} \circ \pi )^{*}\Omega_{X}$.
We can define $\pi_{X}$ with respect to $\Omega_{X}$ for 
the Calabi-Yau manifold ${X}$ as well.

In order to describe the structure of the intermediate Jacobian of $X$ 
we look at the period map of an abelian threefold, 
following the method in Shioda~\cite{Sh1} (also see Mumford~\cite{M}).
Choose a basis $u_1$, $u_2$ of $H_1(E,{\bf R})$. This determines a 
${\bf C}$-basis $H_1(E,{\bf R}) \otimes_{\bf R} {\bf C}= H_1(E,{\bf C})$.
If $e_i$ for $i=1,2,3$ is the standard basis of ${\bf C}^3$ then
$u_{2i-1}=e_i$ and $u_{2i}=\omega e_i$ for $i=1,2,3$ 
form a basis of $H_1(A,{\bf Z})$ and $A= {\bf C}^3/M$ with $M$ the lattice
generated by $u_1,\ldots,u_6$. The dual basis is denoted by $v_i$.
The basis of $H^1(A,{\bf Z})$ determines a canonical basis 
$v_i \wedge v_j \wedge v_k$ of $H^3(A,{\bf Z})$.
The natural homomorphism
$$
p_{A} : H^{3}(A, {\bf Z}) \longrightarrow H_{dR}^{3}(A) \longrightarrow 
H^{3}(A, O_{A}) \cong {\bf C}.
$$
is an element of $\Hom_{\bf C}(H^{3}(A,{\bf C}), {\bf C})$ and can
be considered as an element of $H^{3}(A, {\bf C})$ and is given by
$$
p_{A} = \sum_{i < j < k}\det(u_{i}, u_{j}, u_{k})v^{i}\wedge v^{j}\wedge v^{k}.
$$
Therefore the image of $p_{A}$ in ${\bf C}$ 
is spanned by the complex numbers $1$, $\omega$, $\omega^{2}$ and
$\omega^{3}$ over $\bf Z$.

\begin{lemma}\label{lm;surj}
Let $X$ and $Y$ be topological manifolds, and let
$\rho : X  \longrightarrow Y$ be a unramified surjective continuous map.
Then, $\rho_{*} : H^{i}(X, {\bf Z})  \longrightarrow H^{i}(Y, {\bf Z})$ is
surjective.
\end{lemma}
\proof{By Poincar\'e duality, it suffices to prove that 
$\rho_{*} : H_{i}(X, {\bf Z})  \longrightarrow H_{i}(Y, {\bf Z})$ is
surjective. With a sufficiently small triangulation of $Y$, 
this follows from the fact
that $\rho$ is a local isomorphism.}
 
We now give the proof of Theorem \ref{intermediatejac}.
Let $S$ be a set of non-free points of the action of
$G$ on $A$. Then the restriction of $\pi$ to  $A\setminus S$
is \'etale on $A/G \setminus \pi(S)$.  Since $S$ is of codimension $3$ in $A$, 
we have the following diagram:
$$
\begin{array}{ccccl}
H^{3}(A, {\bf Z}) &\cong & H^{3}(A\setminus S, {\bf Z}) & 
\stackrel{p_{A}}{\longrightarrow} &
H^{3}(A,O_{A}) \cong {\bf C} \\
\downarrow{\pi_{*}}&  &\downarrow{(\pi\mid_{A\setminus S})_{*}} & & 
\quad \downarrow{\pi_{*}}\\
H^{3}(A/G, {\bf Z}) &\cong &H^{3}(A/G \setminus \pi(S), {\bf Z}) &
\stackrel{p_{A/G}}{\longrightarrow}  & H^{3}(A/G, {{O}}_{A/G}) \cong {\bf C}\\
 \downarrow{\cong} & & & &\quad \downarrow{\nu^{*}} \\
H^{3}(X, {\bf Z})& & & \stackrel{p_{X}}{\longrightarrow} & 
H^{3}(X, {{O}}_{X}) 
\cong {\bf C}.
\end{array}
$$
The vertical arrows on the right hand side give an identification
of $H^{3}(A,O_{A})$ and $H^{3}(X,O_{X})$. Since 
$(\pi\mid_{A\setminus S})_{*}$ is surjective by Lemma~\ref{lm;surj},
$\pi_{*}$ is also surjective and  ${\rm Im}\,  p_{A} = 
{\rm Im} \, p_{X}$ in ${\bf C}$. 
Now ${\rm Im} \, p_{X}$ is a lattice in ${\bf C}$, and ${\rm Im}\, p_{A}$ 
is a lattice in ${\bf C}$ as well. We know that 
${\rm Im} \, p_{A}$ is generated by $1$, $\omega$, $\omega^{2}$ 
and $\omega^{3}$ and thus 
$\omega$ is a quadratic number and the intermediate Jacobian
has complex multiplication
by ${\bf Q}(\omega)$.
Hence the intermediate Jacobian ${\bf C}/ {\rm Im}\,  p_{X}$ of $X$ 
is isogenous to $E$.
If $\omega$ is a quadratic integer, then we have 
${\rm Im} \, p_{X} ={\rm Im} p_{A} = {\bf Z} + {\bf Z}\omega$, and 
the intermediate Jacobian ${\bf C}/ {\rm Im}\,  p_{X}$ of $X$ 
is isomorphic to $E$.
\end{section}

\begin{section}{Questions}

We close with two natural basic questions that suggest themselves.

Is there a function $f(n)$ such that a Calabi-Yau variety
in characteristic $p>0$ of dimension $n$
lifts to characteristic $0$ if $p> f(n)$?
Note that Hirokado constructed a non-liftable Calabi-Yau threefold
in characteristic $3$, see \cite{Hir}.

Can a Calabi-Yau variety of dimension $3$ in positive characteristic
have non-zero regular $1$-forms or regular $2$-forms? 
\begin{section}{Acknowledgement}
This research was made possible by a JSPS-NWO grant. The second
author would like to thank the University of Amsterdam for hospitality
and the first author would like to thank Prof.\ Ueno for inviting him
to Kyoto, where this paper was finished.
\end{section}

\end{section}

\end{document}